\documentclass[11pt]{amsart}

\usepackage{amsmath}
\usepackage{amssymb}
\usepackage{amsthm}
\usepackage{enumerate}
\usepackage{amsbsy}
\usepackage{amsfonts}

\usepackage[bitstream-charter]{mathdesign}
\usepackage[T1]{fontenc}

\usepackage{amsthm}
\newtheorem{thm}{Theorem}

\newtheorem{lem}{Lemma}

\theoremstyle{definition}
\newtheorem{defn}{Definition}
\theoremstyle{remark}
\newtheorem{rem}{Remark}

\newcommand{\bN}{\mathbb{N}}
\newcommand{\bR}{\mathbb{R}}
\newcommand{\bU}{\mathbf{U}}
\newcommand{\bB}{\mathbf{B}}
\newcommand{\bP}{\mathrm{P}}
\newcommand{\cC}{\mathcal{C}}
\newcommand{\cS}{\mathcal{S}}
\newcommand{\set}[1]{\{#1\}}

\begin{document}
\baselineskip=18pt

\title{The Minkowski dimension of interior singular points in the incompressible Navier--Stokes equations}

\author{Youngwoo Koh \& Minsuk Yang}

\address{Youngwoo Koh: School of Mathematics, Korea Institute for Advanced Study, 85 Hoegiro Dongdaemungu, Seoul, Republic of Korea}
\address{Minsuk Yang: School of Mathematics, Korea Institute for Advanced Study, 85 Hoegiro Dongdaemungu, Seoul, Republic of Korea}
\email{yangm@kias.re.kr}

\begin{abstract}
We study the possible interior singular points of suitable weak solutions to the three dimensional incompressible Navier--Stokes equations.
We present an improved parabolic upper Minkowski dimension of the possible singular set, which is bounded by $95/63$.
The result also continue to hold for the three dimensional incompressible magnetohydrodynamic equations without any difficulty.\\

\end{abstract}

\maketitle

\section{Introduction}
\label{S1}

We study the possible interior singular points of suitable weak solutions to the three dimensional incompressible Navier--Stokes equations with unit viscosity and zero external force
\begin{align*}
(\partial_t - \Delta) \bU + (\bU\cdot\nabla) \bU + \nabla \bP &= 0 \\
\nabla\cdot\bU &= 0
\end{align*}
in $\Omega_T=\Omega \times (0,T) \subset \bR^3 \times \bR$.
We assume that the initial data $\bU_0$ is sufficiently regular.
The state variables $\bU$ and $\bP$ denote the velocity field of the fluid and its pressure.
In this paper, we only treat local interior regularity theory.

The study of regularity problem for the Navier--Stokes equations has long history and huge literature.
We mention only a few of them.
Scheffer \cite{MR0454426,MR0510154} introduced partial regularity for the Navier--Stokes equations.
Caffarelli, Kohn and Nirenberg \cite{MR673830} further strengthened Scheffer's results.
Lin \cite{MR1488514} gave a new short proof by an indirect argument.
Ladyzhenskaya and Seregin \cite{MR1738171} gave a clear presentation of the H\"older regularity. 
Choe and Lewis \cite{MR1780481} studied the singular set by using a generalized Hausdorff measure.
Gustafson, Kang and Tsai \cite{MR2308753} unified the several known criteria.

We would like to discuss about the two important criteria.
For $z=(x,t) \in \Omega_T$ we denote balls and parabolic cylinders by
\begin{align*}
B(x,r) &= \set{y\in\bR^3:|y-x|<r}, \\
Q(z,r) &= B(x,r) \times (t-r^2,t).
\end{align*}
The first one is as follows.
There exists an absolute positive constant $\epsilon_1$ such that $z$ is a regular point if for some $R>0$
\[
R^{-2} \iint_{Q(z,R)} |\bU|^3 + |\bP|^{3/2} dyds < \epsilon_1.
\]
The second one is as follows.
There exists an absolute positive constant $\epsilon_2$ such that $z$ is a regular point if
\[
\limsup_{r\to0} r^{-1} \int_{Q(z,r)} |\nabla \bU|^2 dyds < \epsilon_2.
\]
The contrapositives of the regularity criteria provide us some information about the possible singularities.
Although the proof of the latter criterion is actually depends on the first criterion, it has an important application that the one dimensional parabolic Hausdorff measure of the possible singular points is zero. 
We say that a point $z \in \Omega_T$ is a singular point of a suitable weak solution $(\bU,\bP)$ if $\bU$ is not locally H\"older continuous in any parabolic neighbourhood $Q(z,r)$.
We denote by $\cS$ the set of all singular points.

There are many different notions measuring lower dimensional sets like $\cS$. 
Another important one is the Minkowski (box-counting) dimension.
When we count the number of uniform boxes to cover the possible singularities, the former criterion gives better information.
If we try to use the second one, then the size of the local cover depends on the point so that we have trouble making a uniform cover.
When one investigate the Minkowski dimension of the singular set, one should use the first kind ``for some $R$'' criterion.
But, it is difficult to lower the Minkowski dimension under $5/3$ due to the scaling structure of the Navier--Stokes equations.
Thus, the natural strategy is combining the different scaled functionals to lower the power of ``$R$''. 

Our objective of this paper is to give an improved bound of the upper Minkowski dimension of the singular set.
We first present Theorem \ref{T1} about the regularities.
And then we use its contraposition to prove Theorem \ref{T2} about the singularities.
We shall give the necessary preliminaries in the next two sections.
Here are the statement of our main results.

\begin{thm}
\label{T1}
For each $\gamma < 10/63$ there exist positive numbers $\varepsilon<1$ and $\rho_0<1$ such that the point $z$ is regular if for some $0 < \rho < \rho_0$
\begin{equation}
\int_{Q(z,\rho)} |\nabla \bU|^2 + |\bU|^{10/3} + |\bP|^{5/3} dxdt
< \rho^{5/3-\gamma} \varepsilon.
\end{equation}
\end{thm}

The contraposition of this regularity criterion yields the following theorem about the possible singularities.

\begin{thm}
\label{T2}
The parabolic upper Minkowski dimension of the set $\cS$ of the possible singular points for the Cauchy problem is bounded by $95/63$.
\end{thm}

\begin{rem}
Note that the bound $95/63$ is better than the previous results $45/29$ in \cite{MR2967124} and $135/82$ in \cite{Ku}.
In \cite{MR2864801} there are several criteria in terms of velocity, velocity gradient, and vorticity, and those criteria can give the bound $5/3$.
\end{rem}

\begin{rem}
Suppose that the existence of suitable weak solutions is established for some bounded domain $\Omega_T = \Omega \times (0,T) \subset \bR^3\times\bR$.
Then for any subdomain $K$ compactly embedded in $\Omega_T$ we carry out the same argument in this paper to obtain the same bound $95/63$ for $K \cap \cS$.
One can cover the compact set $K$ by using a `finite' number of cubes by utilizing Whitney's decomposition.
However, we should handle singular points near the boundary in a different way.
For simplicity we pretend the domain to be $\Omega_T = \bR^4$ so that 
we may restrict the size of various balls in the proof simply by 1.
\end{rem}

\begin{rem}
The existence of a suitable weak solution was proved for the Cauchy problem in \cite{MR673830} and for bounded $C^2$ domain in Section 2.5 of \cite{FKS}.
The existence of a suitable weak solution is closely related to the geometry of the domain.
For the study of more general bounded or unbounded domains, we refer the reader to \cite{FKS} and the references therein.
\end{rem}

\section{Fractal dimensions}
\label{S2}

There are several different ways to measure the amount of sparse sets considering their complex geometric distribution.
Two of the most popular concepts are the Hausdorff dimension and the Minkowski dimension.
We recall here the parabolic versions of the definitions.

\begin{defn}[The parabolic Hausdorff dimension]
For fixed $\delta>0$ and $S\subset\bR^3\times\bR$, let $\cC(S,\delta)$ be the family of all coverings of parabolic cylinders $\set{Q(z_j,r_j)}$ that covers $S$ with $0<r_j \le \delta$.
The $\alpha$ dimensional parabolic Hausdorff measure is defined as
\[
H^\alpha(S) = \lim_{\delta\to0} \inf_{\cC(E,\delta)} \sum_j r_j^\alpha.
\]
The parabolic Hausdorff dimension of the set $S$ is defined as
\[
\dim_H(S) = \inf\set{\alpha : H^\alpha(S)=0}.
\]
\end{defn}

\begin{defn}[The parabolic Minkowski dimension]
Let $N(S;r)$ denote the minimum number of parabolic cylinders $Q(z,r)$ required to cover the set $S$.
Then the parabolic upper Minkowski dimension of the set $S$ is defined as
\begin{equation}
\overline{\dim}_M(S) = \limsup_{r\to0} \frac{\log N(S;r)}{-\log r}
\end{equation}
and the parabolic lower Minkowski dimension of the set $S$ is defined as
\[
\underline{\dim}_M(S) = \liminf_{r\to0} \frac{\log N(S;r)}{-\log r}.
\]
If the limit exists, then it is called the parabolic Minkowski dimension of the set $S$.
\end{defn}

In general, different fractal dimensions of the same set $S$ are not equivalent and they reflect the different geometric structures of the set.
The Minkowski dimension is strongly related to the Hausdorff dimension and a good control of the upper Minkowski dimension has a stronger implication.
Indeed, from the definition it is easy to see that 
\[
\dim_H(S) \le \underline{\dim}_M(S) \le \overline{\dim}_M(S),
\]
but both inequalities may be strict.
For instance, the set $\{1/n: n \in \bN\}$ has the Hausdorff dimension zero but the Minkowski dimension $1/2$.
We refer the reader Falconer's book \cite{MR3236784} for the introduction of the fractal geometry.

\section{Preliminary lemmas}
\label{S3}

We first recall the definition of suitable weak solutions.

\begin{defn}[Suitable weak solutions]
Let $\Omega \subset \bR^3$ and $T>0$.
A pair $(\bU,\bP)$ is called a suitable weak solution in $\Omega_T = \Omega \times (-T,0)$ if the following three conditions are satisfied:
\begin{enumerate}
\item
$\bU \in L^{\infty}(-T,0;L^2(\Omega_T)) \cap L^2(-T,0;H_0^1(\Omega_T))$ and $\bP \in L^{3/2}(\Omega_T)$.
\item
There exists a distribution $\bP$ such that $(\bU,\bP)$ solves the incompressible Navier--Stokes equations in the sense of distributions.
\item
$(\bU,\bP)$ satisfies the local energy inequality: for almost all $t \in (-T,0)$ and for every non-negative $\phi \in C_c^\infty(\bR^3\times\bR)$ vanishing in a neighborhood of the parabolic boundary of $\Omega_T$
\begin{equation}
\begin{split}
\label{E31}
&\int |\bU(\cdot,t)|^2 \phi(\cdot,t) dx + 2 \int_0^t \int_\Omega  |\nabla \bU|^2 \phi dxdt \\
&\le \int_0^t \int_\Omega  |\bU|^2 (\partial_t \phi + \Delta \phi) + (|\bU|^2 + 2\bP) \bU \cdot \nabla \phi dxdt.
\end{split}
\end{equation}
\end{enumerate}
\end{defn}

For notational convenience we shall use the following scaled funcitonals.

\begin{defn}[Scaled functionals]
Let
\begin{align*}
A(z,r) &= r^{-1} \sup_{|t-s|<r^2} \int_{B(x,r)} |\bU(y,s)|^2 dy \\
E(z,r) &= r^{-1} \int_{Q(z,r)} |\nabla \bU(y,s)|^2 dyds \\
C(z,r) &= r^{-2} \int_{Q(z,r)} |\bU(y,s)|^3 dyds \\
D(z,r) &= r^{-2} \int_{Q(z,r)} |\bP(y,s)|^{3/2} dyds.
\end{align*}
\end{defn}

\begin{rem}
The point $z$ in the scale functionals can be suppressed when it is a fixed reference point and understood obviously in the context.
\end{rem}

We end this section by presenting fundamental inequalities between the scaled functionals and a critical regularity criterion, which can be found in many papers concerning partial regularity of the Navier--Stokes equations. 
Interpolation inequalities and pressure inequalities play important roles to complete iteration schemes.
Although there are many variations, interpolation inequalities are basically based on a simple $L^p$ interpolation and the Poincar\'e--Sobolev inequality.
Pressure inequalities depend on decompositions of a localized pressure.
A basic way to get a pressure inequality is to split the localized pressure into the sum of a singular integral and a harmonic function. 
Explicit decomposition of the localized pressure and the idea of interpolation inequality were given in Section 2 and Section 3 of \cite{MR673830}.
For the proof of the following versions we refer the reader to Lemma 2 in \cite{MR3322370} and Lemma 3.4 in \cite{MR2308753}, respectively.

\begin{lem}[Interpolation inequality]
\label{L31}
For any $0<r<1$ and $0<\theta<1$
\[C(\theta r) \le K_1 \theta^{-3/2} A(r)^{3/4} E(r)^{3/4}
+ K_1 \theta^3 A(r)^{3/2}\]
where the constant $K_1>1$ is absolute.
\end{lem}

\begin{lem}[Pressure inequality]
\label{L32}
For any $0<r<1$ and $0<\theta<1/2$
\[D(\theta r) \le K_2 \theta D(r) + K_2 \theta^{-2} C(r)\]
where the constant $K_2>1$ is absolute.
\end{lem}

\begin{lem}[Regularity criterion]
\label{L33}
There exists a positive constant $\zeta$ such that an interior point $z$ is a regular point if for some positive number $R$
\[D(z,R) + C(z,R) \le \zeta.\]
\end{lem}

\begin{rem}
We impose the restriction $0<r<1$ in Lemma \ref{L31} and \ref{L32} for convenience.
\end{rem}

\section{Proof of Theorem \ref{T1}}
\label{S4}

In this section we fix $z$ and suppress it.
Fix a positive number $\gamma < \frac{10}{63}$.
We shall introduce several parameters that depend only on $\gamma$.
First we choose a natural number $N$ such that 
\begin{equation}
\label{E41}
\beta := \frac{1}{6N} < \frac{7}{15} \Big(\frac{10}{63} - \gamma\Big).
\end{equation}
Next, we choose positive numbers $\rho_0$ and $\epsilon$ satisfying $\rho_0^\beta < 1/2$
and 
\begin{equation}
\label{E42}
\varepsilon < \min\left\{\left(\frac{\zeta}{K_2^N K_3+4 K_2^N K_1 K_3^{3/2}}\right)^{10/9},1\right\}
\end{equation}
where $K_1$, $K_2$ and $\zeta$ are the constants in Lemma \ref{L31}, \ref{L32} and \ref{L33} and $K_3=40 (64\pi)^{2/5}$.

\begin{rem}
We note that the parameter $N$ will represent the number of iterations.
So, if $\gamma$ was very close to $10/63$, then $\beta$ should be very small and hence $N$ should be very large. 
It means that we need to perform $N$ iterations.
The particular choice $N\beta=1/6$ looks a little bit mysterious at this moment.
It is an almost optimal choice in our argument based on the inequalities in Lemma \ref{L31} and \ref{L32}.
\end{rem}

Now, we are ready to prove the theorem under this settings.
We divide the proof into a few steps.

\textbf{Step 1)}
Suppose that for some fixed $\rho < \rho_0$ 
\begin{equation}
\label{E43}
\int_{Q(z,2\rho)} |\nabla \bU|^2 + |\bU|^{10/3} + |\bP|^{5/3} 
\le (2\rho)^{5/3-\gamma} \varepsilon.
\end{equation}
We claim that 
\begin{equation}
\label{E44}
\begin{split}
A(\rho) &= \rho^{-1} \sup_{|t-s|<\rho^2} \int_{B(\rho)} |\bU|^2 
\le K_3 \rho^{-9\gamma/10} \varepsilon^{3/5}, \\
E(\rho) &= \rho^{-1} \int_{Q(\rho)} |\nabla \bU|^2 
\le 2^{5/3} \rho^{2/3-\gamma} \varepsilon.
\end{split}
\end{equation}

Indeed, the estimate of $E(\rho)$ follows immediately from its definition and \eqref{E43}.
On the other hand, for the estimate $A(\rho)$ we use the local energy inequality with a cutoff function $\phi$, which is smooth, supported in $Q(z,2\rho)$, $\phi \geq 1$ in $Q(z,\rho)$ and 
\[|\partial_t\phi|+|\nabla^2\phi|+|\nabla\phi|^2 \le 10\rho^{-2} \quad \text{ in } Q(z,2\rho).\]
From the local energy inequality \eqref{E31}, like Lemma 5.2 in \cite{MR1738171}, we use the H\"older inequality, the assumption \eqref{E43}, $\varepsilon<1$ and $\gamma>0$ to obtain that 
\begin{align*}
A(\rho) 
&\le 10 \rho^{-3} \int_{Q(2\rho)} |\bU|^2
+ 10 \rho^{-2} \int_{Q(2\rho)} |\bU|^3
+ 10 \rho^{-2} \int_{Q(2\rho)} |\bP\bU| \\
&\le 10 (64\pi)^{2/5} \rho^{-1} \Big(\int_{Q(2\rho)} |\bU|^{10/3}\Big)^{3/5}
+ 10 (64\pi)^{1/10} \rho^{-3/2} \Big(\int_{Q(2\rho)} |\bU|^{10/3}\Big)^{9/10} \\
&\quad + 10 (64\pi)^{1/10} \rho^{-3/2} \Big(\int_{Q(2\rho)} |\bU|^{10/3}\Big)^{3/10} \Big(\int_{Q(2\rho)} |\bP|^{5/3}\Big)^{3/5} \\
&\le 20 (64\pi)^{2/5} 
\Big(\rho^{-3\gamma/5} \varepsilon^{3/5}
+ \rho^{-9\gamma/10} \varepsilon^{9/10}\Big) \\
&\le 40 (64\pi)^{2/5} \rho^{-9\gamma/10} \varepsilon^{3/5}. 
\end{align*}
Since $K_3=40 (64\pi)^{2/5}$, this proves the claim \eqref{E44}.

\textbf{Step 2)}
For $j=0,1,2,\dots,N$ we define 
\begin{equation}
\label{E45}
R_j := \rho^{\alpha+j\beta} \quad \text{and} \quad \theta := \rho^\beta < \rho_0^\beta < 1/2
\end{equation}
where 
\begin{equation}
\label{E46}
\alpha := \frac{550-9\gamma}{540} > 1.
\end{equation}
We note that $R_N < R_{N-1} < \dots < R_0 = \rho^\alpha < \rho_0^\alpha < \rho_0 < 1$.
By the pressure inequality (Lemma \ref{L32}), we have for $j=0,1,\dots,N-1$
\[D(R_{j+1}) \le K_2 \theta D(R_j) + K_2 \theta^{-2} C(R_j).\]
Iterating the pressure inequality (Lemma \ref{L32}), we obtain that  
\[
D(R_N) + C(R_N) \le (K_2\theta)^N D(R_0) + \sum_{j=0}^{N} (K_2\theta)^{N-j} \theta^{-3} C(R_j),
\]
which is obvious from an induction argument on $N$.
It suffices to show that 
\begin{equation}
\label{E47}
I := (K_2\theta)^N D(R_0) \le K_2^N K_3 \varepsilon^{9/10}
\end{equation}
and 
\begin{equation}
\label{E48}
II := \sum_{j=0}^{N} (K_2\theta)^{N-j} \theta^{-3} C(R_j) 
\le 4 K_2^N K_1 K_3^{3/2} \varepsilon^{9/10}.
\end{equation}
We prove them in the next steps, respectively.
Suppose, for the moment, the estimates of $I$ and $II$ are obtained.
From \eqref{E42} we have
\[D(R_N) + C(R_N) \le \big(K_2^N K_3+4 K_2^N K_1 K_3^{3/2}\big) \varepsilon^{9/10} < \zeta.\]
Therefore, $z$ is a regular point due to Lemma \ref{L33}. 

\textbf{Step 3)}
We recall $R_0 = \rho^\alpha$ and $\theta = \rho^\beta$.
Using the H\"older inequality and the assumption \eqref{E43} we have 
\begin{align*}
I 
&= (K_2\theta)^N D(R_0) \\
&\le K_2^N (64\pi)^{1/10} \theta^N R_0^{-3/2} \Big(\int_{Q(R_0)} |\bP|^{5/3}\Big)^{9/10} \\
&\le K_2^N 4(64\pi)^{1/10} \theta^N R_0^{-3/2} \Big(\varepsilon \rho^{5/3-\gamma}\Big)^{9/10} \\
&\le K_2^N K_3 \rho^{3/2-3\alpha/2+N\beta-9\gamma/10} \varepsilon^{9/10}.
\end{align*}
From \eqref{E41}, \eqref{E46}, and $\gamma < \frac{10}{63}$, we see that 
\[3/2-3\alpha/2+N\beta-9\gamma/10 = \frac{10-63\gamma}{72} > 0\]
and hence $\rho^{3/2-3\alpha/2+N\beta-9\gamma/10} \le \rho_0^{\frac{10-63\gamma}{72}} \le 1$.
This proves the estimate \eqref{E47}.

\textbf{Step 4)}
Using the interpolation inequality (Lemma \ref{L31}) and the estimates in \eqref{E44}, we have 
\begin{align*}
C(R_j) 
&\le K_1 \Big(\frac{\rho}{R_j}\Big)^{3/2} A(\rho)^{3/4} E(\rho)^{3/4}
+ K_1 \Big(\frac{R_j}{\rho}\Big)^{3} A(\rho)^{3/2} \\
&\le K_1 \Big(\frac{\rho}{\rho^{\alpha+j\beta}}\Big)^{3/2} 
\Big(K_3 \varepsilon^{3/5} \rho^{-9\gamma/10}\Big)^{3/4} 
\Big(K_3 \varepsilon \rho^{2/3-\gamma}\Big)^{3/4} \\
&\quad + K_1 \Big(\frac{\rho^{\alpha+j\beta}}{\rho}\Big)^{3}
\Big(K_3 \varepsilon^{3/5} \rho^{-9\gamma/10}\Big)^{3/2} \\
&\le K_1 K_3^{3/2} \varepsilon^{9/10} 
\Big(\rho^{2-3\alpha/2-3\beta j/2-57\gamma/40}
+ \rho^{-3+3\alpha+3\beta j-27\gamma/20}\Big).
\end{align*}
Hence, using the above inequality, summing geometric series, and substituting $\alpha$ defined in \eqref{E46}, we obtain that 
\begin{align*}
II
&\le K_2^N \sum_{j=0}^{N} \theta^{N-j} \theta^{-3} C(R_j) \\
&\le K_2^N K_1 K_3^{3/2} \varepsilon^{9/10} \sum_{j=0}^{N} \theta^{N-j} \theta^{-3} 
\Big(\rho^{2-3\alpha/2-3\beta j/2-57\gamma/40}
+ \rho^{-3+3\alpha+3\beta j-27\gamma/20}\Big) \\
&\le K_2^N K_1 K_3^{3/2} \varepsilon^{9/10} \rho^{N\beta-3\beta}
\Big(\rho^{2-3\alpha/2-57\gamma/40} \sum_{j=0}^{N} \rho^{-5\beta j/2}
+ \rho^{-3+3\alpha-27\gamma/20} \sum_{j=0}^{N} \rho^{2\beta j}\Big) \\
&\le 2 K_2^N K_1 K_3^{3/2} \varepsilon^{9/10} \rho^{N\beta-3\beta}
\Big(\rho^{2-3\alpha/2-57\gamma/40-5N\beta/2} 
+ \rho^{-3+3\alpha-27\gamma/20} \Big) \\
&= 4 K_2^N K_1 K_3^{3/2} \varepsilon^{9/10} 
\rho^{2/9-3\beta-7\gamma/5}. 
\end{align*}
From \eqref{E41} we see that 
\[2/9-3\beta-7\gamma/5 = \frac{10-63\gamma}{45} -3\beta > 0\]
and hence $\rho^{2/9-3\beta-7\gamma/5} \le \rho_0^{(1-63\gamma)/45-3\beta} \le 1$.
Therefore we get the estimate \eqref{E48}.

This completes the proof of Theorem \ref{T1}.

\section{Proof of Theorem \ref{T2}}
\label{S5}

For the Cauchy problem it suffices to estimate the Minkowski dimension of the set 
\[
\cS \cap [0,1]^4
\]
since all the estimates have the translation invariant bounds.
Theorem \ref{T1} implies that if $z$ is an interior singular point, 
then for all $r < \rho_0$
\[
\varepsilon r^{5/3-\gamma}
\le \int_{Q(z,r)} |\nabla \bU|^2 + |\bU|^{10/3} + |\bP|^{5/3}.
\]
Now, fix $5r < \rho_0$ and consider the covering $\set{Q(z,r) : z \in \cS}$.
By the Vitali covering lemma, there is a finite disjoint sub-family
\[
\set{Q(z_j,r) : j=1,2,\dots,M}
\]
such that $\cS \subset \bigcup Q(z_j,5r)$.
Summing the inequality above at $z_j$ for $j=1,2,\dots,M$ yields
\begin{align*}
M \varepsilon  r^{5/3-\gamma }
&\le \sum_{i=1}^{M} \int_{Q(z_j,r)} |\nabla \bU|^2 + |\bU|^{10/3} + |\bP|^{5/3} \\
&\le \int_{\Omega_T} |\nabla \bU|^2 + |\bU|^{10/3} + |\bP|^{5/3} =: K_4 < \infty.
\end{align*}
Let $N(\cS \cap [0,1]^4;r)$ denote the minimum number of parabolic cylinders $Q(z,r)$ required to cover the set $\cS \cap [0,1]^4$.
Then 
\[
N(\cS \cap [0,1]^4;r) \le M \le K_4 \varepsilon^{-1} r^{-5/3+\gamma}
\]
and hence 
\[
\limsup_{r\to0} \frac{\log N(\cS \cap [0,1]^4;r)}{-\log r} \le 5/3-\gamma.
\]
Since $\gamma$ can be arbitrarily close to $10/63$, this completes the proof of Theorem \ref{T2}.

\begin{rem}
Suppose $\Omega_T$ is a bounded domain and $K$ is a compactly embedded subdomain.
Let $\set{C_1, C_2, \dots, C_M}$ be a set of finite cubes that covers $K$ and denote by $N(C_k;r)$ the minimum number of parabolic cylinders $Q(z,r)$ required to cover the set $\cS \cap C_k$.
Then it is easy to see that
\[\log N(K;r) \le \log \bigg(\sum_{k=1}^M N(C_k;r)\bigg) \le \log M + \max_{1 \le k \le M} \set{\log N(C_k;r)}\]
and therefore
\[
\limsup_{r\to0} \frac{\log N(K;r)}{-\log r} \le 5/3-\gamma.
\]
\end{rem}

\section{Magnetohydrodynamic equations}
\label{S6}

This section is a kind of remark.
We refer the reader to Davidson's monograph \cite{MR1825486} for the background material of magnetohydrodynamics.
The three dimensional incompressible magnetohydrodynamic equations described by 
\begin{align*}
&\partial_t \bU - \Delta \bU + (\bU\cdot\nabla) \bU - (\bB\cdot\nabla)\bB + \nabla \bP = 0 \\
&\partial_t \bB - \Delta \bB + (\bU\cdot\nabla) \bB - (\bB\cdot\nabla)\bU = 0 \\
&\nabla\cdot\bU = 0, \quad \nabla\cdot\bB = 0, 
\end{align*}
has a suitable weak solution satisfying the local energy inequality (see e.g. Definition 2.1 in \cite{MR2561280})
\begin{align*}
&\int_\Omega (|\bU|^2+|\bB|^2) \phi dx + 2 \int_0^t \int_\Omega (|\nabla \bU|^2+|\nabla \bB|^2) \phi dxdt \\
&\le \int_0^t \int_\Omega (|\bU|^2+|\bB|^2) (\phi_s+\Delta\phi) dxdt \\
&\quad + \int_0^t \int_\Omega (|\bU|^2+|\bB|^2+2\bP) \bU \cdot \nabla \phi - 2(\bB \cdot \bU) (\bB \cdot \nabla \phi) dxdt.
\end{align*}
Moreover, Lemma \ref{L31}, \ref{L32}, and \ref{L33} continue to hold with the replaced scaled functionals
\begin{align*}
A(z,r) &= r^{-1} \sup_{|t-s|<r^2} \int_{B(x,r)} |\bU(y,s)|^2+|\bB(y,s)|^2 dy \\
E(z,r) &= r^{-1} \int_{Q(z,r)} |\nabla \bU(y,s)|^2+|\nabla \bB(y,s)|^2 dyds \\
C(z,r) &= r^{-2} \int_{Q(z,r)} |\bU(y,s)|^3+|\bB(y,s)|^3 dyds \\
D(z,r) &= r^{-2} \int_{Q(z,r)} |\bP(y,s)|^{3/2} dyds.
\end{align*}
We refer the reader for the interpolation inequality of MHD to Lemma 2 in \cite{MR3322370} and for the pressure inequality of MHD to Lemma 3.3 in \cite{MR2561280}.
They are the essential ingredients used in the proofs.
By the same way in Section \ref{S4} one can easily prove the following theorem.

\begin{thm}
For each $\gamma < 10/63$ there exist positive numbers $\varepsilon<1$ and $\rho_0<1$ such that the point $z$ is regular if for some $0 < \rho < \rho_0$
\[
\int_{Q(z,\rho)} \big(|\nabla \bU|^2 + |\nabla \bB|^2\big) + \big(|\bU|^{10/3} + |\bB|^{10/3}\big) + |\bP|^{5/3} dxdt < \rho^{5/3-\gamma} \varepsilon .
\]
\end{thm}

This theorem has a direct application to the parabolic upper Minkowski dimension of suitable weak solution to the MHD equations and we omit the tedious repetitions.
One may find the necessary information, for example, in the papers \cite{MR2165089}, \cite{MR2561280}, \cite{MR3097270}, \cite{MR3322370} and the references therein.

\end{document}